\documentclass[11pt]{article}

\usepackage{amsmath}
\usepackage{amsthm}
\usepackage{amssymb}
\usepackage{latexsym}
\usepackage{pstricks}
\usepackage{graphicx, pst-plot, pst-node, pst-text, pst-tree}
\usepackage{titlefoot}
\usepackage{enumerate}
\usepackage{titlesec}
\usepackage{mathabx}
\usepackage[small,it]{caption}
\usepackage{color}
\definecolor{lightgray}{rgb}{0.8, 0.8, 0.8}
\definecolor{darkgray}{rgb}{0.7, 0.7, 0.7}
\definecolor{darkblue}{rgb}{0, 0, .4}

\usepackage[bookmarks,breaklinks]{hyperref}
\hypersetup{
        colorlinks=true,
        linkcolor=darkblue,
        anchorcolor=darkblue,
        citecolor=darkblue,
        pagecolor=darkblue,
        urlcolor=darkblue,
        pdfpagemode=UseThumbs,
        pdftitle={A sharp bound for the reconstruction of partitions},
        pdfsubject={Combinatorics},
        pdfauthor={Vincent Vatter},
}

\newtheorem{theorem}{Theorem}

\newtheoremstyle{example}{\topsep}{\topsep}%
     {}
     {}
     {\bfseries}
     {.}
     {.5em}
     {\thmname{#1}\thmnumber{ #2}}
\theoremstyle{example}

\newtheoremstyle{negexample}{\topsep}{\topsep}%
     {}
     {}
     {\bfseries}
     {.}
     {.5em}
     {\thmname{#1}\thmnumber{ #2}}
\theoremstyle{negexample}
\newtheorem{negexample}[theorem]{Negative Example}

\newcounter{todocounter}

\long\def\symbolfootnote[#1]#2{\begingroup%
\def\thefootnote{\fnsymbol{footnote}}\footnote[#1]{#2}\endgroup}

\makeatletter
\newcommand{\Rm}[1]{\expandafter\@slowromancap\romannumeral #1@}
\newfont{\footsc}{cmcsc10 at 8truept}
\newfont{\footbf}{cmbx10 at 8truept}
\newfont{\footrm}{cmr10 at 10truept}
\renewenvironment{abstract}%
                {
                  \begin{list}{}%
                     {\setlength{\rightmargin}{1in}%
                      \setlength{\leftmargin}{1in}}%
                   \item[]\ignorespaces\begin{small}}%
                 {\end{small}\unskip\end{list}}

\datefoot{\today}
\amssubj{05A17, 06A07}

\newpagestyle{main}[\small]{
        \headrule
        \sethead[\usepage][][]
        {\sc A Sharp Bound for the Reconstruction of Partitions}{}{\usepage}}

\title{\sc{A Sharp Bound for the Reconstruction of Partitions}}
\author{\sc{Vincent Vatter}\\
\small Department of Mathematics\\[-1pt]
\small Dartmouth College\\[-1pt]
\small Hanover, NH 03755\\[-10pt]}

\titleformat{\section}
        {\large\sc}
        {\thesection.}{1em}{}   

\date{}

\begin{document}
\maketitle

\pagestyle{main}

\begin{abstract}
Answering a question of Cameron, Pretzel and Siemons proved that every integer partition of $n\ge 2(k+3)(k+1)$ can be reconstructed from its set of $k$-deletions.  We describe a new reconstruction algorithm that lowers this bound to $n\ge k^2+2k$ and present examples showing that this bound is best possible.
\end{abstract}

Analogues and variations of Ulam's notorious graph reconstruction conjecture have been studied for a variety of combinatorial objects, for instance words (see Sch\"utzenberger and Simon~\cite[Theorem 6.2.16]{lothaire:combinatorics-o:}), permutations (see Raykova~\cite{raykova:permutation-rec:} and Smith~\cite{smith:permutation-rec:}), and compositions (see Vatter~\cite{vatter:reconstructing-:}), to name a few.  

In answer to Cameron's query~\cite{cameron:stories-from-th:} about the partition context, Pretzel and Siemons~\cite{pretzel:reconstruction-:} proved that every partition of $n\ge 2(k+3)(k+1)$ can be reconstructed from its set of $k$-deletions.  Herein we describe a new reconstruction algorithm%
%
\ that lowers this bound, establishing the following result, which Negative Example~\ref{exa-part-badpair} shows is best possible.

\begin{theorem}\label{thm-part-reconstruct}
Every partition of $n\ge k^2+2k$ can be reconstructed from its set of $k$-deletions.
\end{theorem}

We begin with notation.  Recall that a {\it partition of $n$\/}, $\lambda=(\lambda_1,\dots,\lambda_\ell)$, is a finite sequence of nonincreasing integers whose sum, which we denote $|\lambda|$, is $n$.  The {\it Ferrers diagram\/} of $\lambda$, which we often identify with $\lambda$, consists of $\ell$ left-justified rows where row $i$ contains $\lambda_i$ cells.  An {\it inner corner\/} in this diagram is a cell whose removal leaves the diagram of a partition, and we refer to all other cells as {\it interior cells\/}.

We write $\mu\le\lambda$ if $\mu_i\le\lambda_i$ for all $i$; another way of stating this is that $\mu\le\lambda$ if and only if $\mu$ is contained in $\lambda$ (here identifying partitions with their diagrams).  If $\mu\le\lambda$, we write $\lambda/\mu$ to denote the set of cells which lie in $\lambda$ but not in $\mu$.  We say that the partition $\mu$ is a {\it $k$-deletion\/} of the partition of $\lambda$ if $\mu\le\lambda$ and $|\lambda/\mu|=k$.

Recall that this order defines a lattice on the set of all finite partitions, known as {\it Young's lattice\/}, and so every pair of partitions has a unique {\it join\/} (or {\it least upper bound\/})
$$
\mu\vee\lambda=(\max\{\mu_1,\lambda_1\},\max\{\mu_2,\lambda_2\},\dots)
$$
and {\it meet\/}
$$
\mu\wedge\lambda=(\min\{\mu_1,\lambda_1\},\min\{\mu_2,\lambda_2\},\dots).
$$

Finally, recall that the {\it conjugate\/} of a partition $\lambda$ is the partition $\lambda'$ obtained by flipping the diagram of $\lambda$ across the NW-SE axis; it follows that $\lambda'_i$ counts the number of entries of $\lambda$ which are at least $i$.

Before proving Theorem~\ref{thm-part-reconstruct} we show that it is best possible:

\begin{negexample}\label{exa-part-badpair}
For $k\ge 1$, consider the two partitions
\begin{eqnarray*}
\mu&=&(\underbrace{k+1,\dots,k+1}_{k},k-1)\mbox{ and}\\
\lambda&=&(\underbrace{k+1,\dots,k+1}_{k-1},k,k).
\end{eqnarray*}
Note that no $k$-deletion of $\mu$ can contain the cell $(k,k+1)$ and that no $k$-deletion of $\lambda$ can contain the cell $(k+1,k)$.  Therefore every $k$-deletion of $\mu$ and of $\lambda$ is actually a $(k-1)$-deletion of
$$
\mu\wedge\lambda=(\underbrace{k+1,\dots,k+1}_{k-1},k,k-1),
$$
so $\mu$ and $\lambda$ cannot be differentiated by their sets of $k$-deletions.
\end{negexample}

We are now ready to prove our main result.

\newenvironment{proof-of-thm-part-reconstruct}{\medskip\noindent {\it Proof of Theorem \ref{thm-part-reconstruct}.\/}}{\qed}
\begin{proof-of-thm-part-reconstruct}
Suppose that we are given a positive integer $k$ and a set $\Delta$ of $k$-deletions of some (unknown) partition $\lambda$ of $n\ge k^2+2k$.  Our goal is to determine $\lambda$ from this information.  We begin by setting $\mu = \bigvee_{\delta\in\Delta} \delta$, noting that we must have $\lambda\ge\mu$.  Hence if $|\mu|=n$ then we have $\lambda=\mu$ and we are immediately done, so we will assume that $|\mu|<n$.

First consider the case where $\mu$ has less than $k$ rows.  Let $r$ denote the bottommost row of $\mu$ which contains at least $k$ cells ($r$ must exist because $\mu$ has less than $k$ rows and $|\mu|\ge k^2+k$).  Thus the $r$th row of $\lambda$ contains at least $k$ cells as well, so there are $k$-deletions of $\lambda$ in which the removed cells all lie in or below row $r$.  Hence the first $r-1$ rows of $\lambda$ and $\mu$ agree.  Now note that $\lambda$ has more than $2k$ cells to the right of column $k$, so there are $k$-deletions of $\lambda$ in which the removed cells all lie to the right of column $k$, and thus the first $k$ columns of $\lambda$ and $\mu$ agree.  This implies that $\lambda$ and $\mu$ agree on all rows below $r$ (since these rows have less than $k$ cells in $\mu$) and so all cells of $\lambda/\mu$ must lie in row $r$, uniquely determining $\lambda$, as desired.  The case where $\mu$ has less than $k$ columns follows by symmetry.

We may now assume that $\mu$ has at least $k$ rows and $k$ columns.  Let $r$ (resp.~$c$) denote the bottommost row (resp.~rightmost column) containing at least $k$ cells.  Both $r$ and $c$ exist because $\mu$ has at least $k$ rows and columns.  Therefore both $\lambda$ and $\mu$ can be divided into three quadrants, \Rm{1}, \Rm{2}, and \Rm{3}, as shown in Figure~\ref{fig-case1-quads}. 

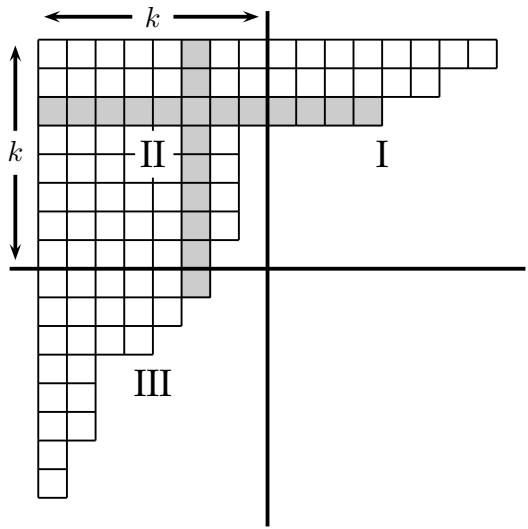
\begin{figure}
\begin{center}
\psset{xunit=0.015in, yunit=0.015in}
\psset{linewidth=0.01in}
\begin{pspicture}(-10,-10)(180,180)
\psframe[linewidth=0,fillstyle=solid,fillcolor=lightgray](0,130)(120,140)
\psframe[linewidth=0,fillstyle=solid,fillcolor=lightgray](50,160)(60,70)
\psline(0,160)(160,160)
\psline(0,150)(160,150)
\psline(0,140)(140,140)
\psline(0,130)(120,130)
\psline(0,120)(70,120)
\psline(0,110)(70,110)
\psline(0,100)(70,100)
\psline(0,90)(70,90)
\psline(0,80)(60,80)
\psline(0,70)(60,70)
\psline(0,60)(50,60)
\psline(0,50)(40,50)
\psline(0,40)(20,40)
\psline(0,30)(20,30)
\psline(0,20)(20,20)
\psline(0,10)(10,10)
\psline(0,0)(10,0)
\psline(0,0)(0,160)
\psline(10,160)(10,0)
\psline(20,160)(20,20)
\psline(30,160)(30,50)
\psline(40,160)(40,50)
\psline(50,160)(50,60)
\psline(60,160)(60,70)
\psline(70,160)(70,90)
\psline(80,160)(80,130)
\psline(90,160)(90,130)
\psline(100,160)(100,130)
\psline(110,160)(110,130)
\psline(120,160)(120,130)
\psline(130,160)(130,140)
\psline(140,160)(140,140)
\psline(150,160)(150,150)
\psline(160,160)(160,150)
\psline[linewidth=0.02in]{<->}(-8,83)(-8,158)
\psframe[linewidth=0,linecolor=white,fillstyle=solid,fillcolor=white](-9,113)(-7,128)
\rput[c](-8,121){$k$}
\psline[linewidth=0.02in]{<->}(2,168)(77,168)
\psframe[linewidth=0,linecolor=white,fillstyle=solid,fillcolor=white](33,167)(47,169)
\rput[c](40,168){$k$}
\psline[linewidth=0.02in](-10,80)(170,80)
\psline[linewidth=0.02in](80,-10)(80,170)
\psframe[linewidth=0,linecolor=white,fillstyle=solid,fillcolor=white](34,112)(47,128)
\rput[c](40,120){\Large\Rm{2}}
\rput[c](120,120){\Large\Rm{1}}
\rput[c](40,40){\Large\Rm{3}}
\end{pspicture}
\end{center}
\caption{An example partition $\mu$ from Case 1 of the proof of Theorem~\ref{thm-part-reconstruct}, divided into three quadrants.  Here $k=8$, and $r$ and $c$ appear shaded.}\label{fig-case1-quads}
\end{figure}

As before, we see that the first $r-1$ rows and $c-1$ columns of $\lambda$ and $\mu$ agree.  We consider three cases based on whether and where $r$ and $c$ intersect.

\newenvironment{proof-case}[1]{\medskip\noindent {\it Case #1.\/}}{}
\begin{proof-case}{1: $r$ and $c$ intersect at an interior cell of $\mu$}
Suppose that $r$ and $c$ intersect at the cell $(i,j)$.  It follows from the maximality of $r$ and $c$ that $i,j<k$, and thus the cell $(k,k)$ does not lie in $\mu$.  Were the cell $(k,k)$ to lie in $\lambda$ then, because $|\lambda|\ge k^2+2k$, $\lambda$ must contain at least $2k$ cells to the right of or below $(k,k)$ and thus $\lambda$ would contain a $k$-deletion with the cell $(k,k)$, a contradiction; thus $\lambda$ also does not contain $(k,k)$.

Hence Quadrant \Rm{2} of $\lambda$ contains less than $k^2$ cells, so $\lambda$ must have more than $k$ cells in quadrant \Rm{1} or \Rm{3}.  Hence there are $k$-deletions of $\lambda$ with more than $k$ cells in quadrant \Rm{1} or \Rm{3}; suppose by symmetry that $\lambda$ and $\mu$ both have more than $k$ cells in quadrant \Rm{1}.

There are then $k$-deletions of $\lambda$ in which the removed cells are all chosen from quadrant \Rm{1}, so $\lambda$ and $\mu$ agree on all cells in quadrants \Rm{2} and \Rm{3}.  This shows that $r$ is also the bottommost row of $\lambda$ with at least $k$ cells, and so $\lambda/\mu$ contains no cells below row $r$ in quadrant \Rm{1}.  As we already know that $\lambda$ and $\mu$ agree on their first $r-1$ rows, we can therefore conclude that all cells of $\lambda/\mu$ lie in row $r$, which allows us to reconstruct $\lambda$ and complete the proof of this case.
\end{proof-case}

\begin{proof-case}{2: $r$ and $c$ intersect at an inner corner of $\mu$}
Then this inner corner must be the rightmost cell of row $r$ and the bottom cell of column $c$.  It follows that $r,c\ge k$.  Because $\lambda$ and $\mu$ agree to the left of column $c$ and above row $r$, all cells of $\lambda/\mu$ must lie below or to the right of $(r,c)$.  However, the cell $(r+1,c+1)$ cannot lie in $\lambda$ because if it did then one could form a $k$-deletion of $\lambda$ by removing only points lying to the right of column $c$, which would leave at least $k$ cells in row $r+1$ and contradict the definition of $r$.  This leaves only two possibilities for $\lambda/\mu$: the cells $(r,c+1)$ and $(r+1,c)$.  However, only one of these cells can be added to $\mu$ to produce a partition; if both could be added then row $r+1$ and column $c+1$ of $\lambda$ would each contain at least $k$ cells, implying the existence of $k$-deletions of $\lambda$ in which each contain at least $k$ cells and thus contradicting the choice of $r$ and $c$.  This case therefore reduces to checking which one of the cells $(r,c+1)$ and $(r+1,c)$ can be added to $\mu$ to produce a partition.
\end{proof-case}

\begin{proof-case}{3: $r$ and $c$ do not intersect}
Suppose that the rightmost cell in row $r$ is $(r,j)$ and the bottommost cell in column $c$ is $(i,c)$.  If $j<c-1$ then because $\lambda$ and $\mu$ agree to the left of column $c$, $\lambda/\mu$ cannot contain any cells in or below row $r$, and we already have that $\lambda$ and $\mu$ agree above row $r$, so we are left with the conclusion that $\lambda=\mu$.  By symmetry we are also done if $i<r-1$, leaving us to consider the case where $i=r-1$ and $j=c-1$.  Again using the fact that $\lambda$ and $\mu$ agree above row $r$ and to the left of column $c$ (and the definitions of $r$ and $c$) we see that the only possibility for $\lambda/\mu$ is $(r,c)$, completing the proof of this case and the theorem.
\end{proof-case}
\end{proof-of-thm-part-reconstruct}

\bigskip\noindent{\bf Acknowledgements.} I would like to thank the referee for several suggestions which improved the transparency of the proof.

\bibliographystyle{acm}
\bibliography{../refs}

\def\cprime{$'$}
\begin{thebibliography}{1}

\bibitem{cameron:stories-from-th:}
{\sc Cameron, P.~J.}
\newblock Stories from the age of reconstruction.
\newblock {\em Congr. Numer. 113\/} (1996), 31--41.

\bibitem{lothaire:combinatorics-o:}
{\sc Lothaire, M.}
\newblock {\em Combinatorics on Words}, vol.~17 of {\em Encyclopedia of
  Mathematics and its Applications}.
\newblock Addison-Wesley Publishing Co., Reading, Mass., 1983.

\bibitem{pretzel:reconstruction-:}
{\sc Pretzel, O., and Siemons, J.}
\newblock Reconstruction of partitions.
\newblock {\em Electron. J. Combin. 11}, 2 (2004--06), Note 5, 6 pp.

\bibitem{raykova:permutation-rec:}
{\sc Raykova, M.}
\newblock Permutation reconstruction from minors.
\newblock {\em Electron. J. Combin. 13\/} (2006), Research paper 66, 14 pp.

\bibitem{smith:permutation-rec:}
{\sc Smith, R.}
\newblock Permutation reconstruction.
\newblock {\em Electron. J. Combin. 13\/} (2006), Note 11, 8 pp.

\bibitem{vatter:reconstructing-:}
{\sc Vatter, V.}
\newblock Reconstructing compositions.
\newblock {\em Discrete Math. 308}, 9 (2008), 1524--1530.

\end{thebibliography}

\end{document}